\documentclass[12pt]{article}

\usepackage[margin=1in]{geometry}

\usepackage{amssymb}
\usepackage{amsthm}
\usepackage{amsmath}
\usepackage{graphicx}
\usepackage{acronym}
\usepackage{amsfonts}
\usepackage{amscd}

\theoremstyle{plain}
\newtheorem{theorem}{Theorem} 
\newtheorem{Example}{Example}
\newtheorem{definition}[theorem]{Definition}
\def\be{\begin{Example}}
\def\ee{\end{Example}}
\def\bt{\begin{theorem}}
\def\et{\end{theorem}\bigskip}
\def\bl{\begin{Lemma}}
\def\el{\end{Lemma}\bigskip}
\def\ep{\end{Proposition}\bigskip}
\def\bp{\begin{Proposition}}
\def\bd{\begin{definition}}
\def\ed{\end{definition}}
\newcommand{\A}{{\cal A}}
\newcommand{\B}{{\cal B}}

\newcommand{\App}{App}

\def\bt{\beta}

\newcommand{\bdes}{\begin{description}}
\newcommand{\edes}{\end{description}}

\newcommand{\bnum}{\begin{enumerate}}
\newcommand{\enum}{\end{enumerate}}

\newcommand{\bit}{\begin{itemize}}
\newcommand{\eit}{\end{itemize}}

\newcommand{\bea}{\begin{eqnarray}}
\newcommand{\eea}{\end{eqnarray}}
\newcommand{\beq}{\begin{equation}}
\newcommand{\eeq}{\end{equation}}

\newcommand{\baray}{\begin{array}}
\newcommand{\earay}{\end{array}}

\newcommand{\bsry}{\begin{subarray}}
\newcommand{\esry}{\end{subarray}}

\newcommand{\bca}{\begin{cases}}
\newcommand{\eca}{\end{cases}}

\newcommand{\bcen}{\begin{center}}
\newcommand{\ecen}{\end{center}}

\newcommand{\bbm}{\begin{bmatrix}}
\newcommand{\ebm}{\end{bmatrix}}

\newcommand{\bmx}{\begin{matrix}}
\newcommand{\emx}{\end{matrix}}

\newcommand{\bpm}{\begin{pmatrix}}
\newcommand{\epm}{\end{pmatrix}}

\newcommand{\btab}{\begin{tabular}}
\newcommand{\etab}{\end{tabular}}

\begin{document}

\title{\bf The Spectral Theory of Tensors \\ (Rough Version)}

\author{
Liqun Qi\thanks{Department of Applied Mathematics, The Hong Kong
Polytechnic University, Hung Hom, Kowloon, Hong Kong. E-mail:
\emph{maqilq@polyu.edu.hk}. His work is supported by the Hong Kong
Research Grant Council.}}

\date{\today}
\maketitle

The spectral theory of tensors is an important part of numerical
multi-linear algebra, or tensor computation \cite{KB, QSW, Va}.

It is possible that the ideas of eigenvalues of tensors had been
raised earlier.   However, it was in 2005, the papers of Lim and Qi
initiated the rapid developments of the spectral theory of tensors.
In 2005, independently, Lim \cite{Li} and Qi \cite{Qi} defined
eigenvalues and eigenvectors of a real symmetric tensor, and
explored their practical application in determining positive
definiteness of an even degree multivariate form. This work extended
the classical concept of eigenvalues of square matrices, forms an
important part of numerical multi-linear algebra, and has found
applications or links with automatic control, statistical data
analysis, optimization, magnetic resonance imaging, solid mechanics,
quantum physics, higher order Markov chains, spectral hypergraph
theory, Finsler geometry, etc, and attracted attention of
mathematicians from different disciplines.   After six years'
developments, we may classify the spectral theory of tensors to 18
research topics.  Before describing these 18 research topics in four
groups, we state the basic definitions and properties of eigenvalues
of tensors.

\medskip

{\bf A. Basic Theory}

An $n$-dimensional homogeneous polynomial form of degree $m$,
$f(x)$, where $x \in \Re^n$, is equivalent to the tensor product of
a {\bf symmetric} $n$-dimensional tensor $\A = \left(a_{i_1, \cdots,
i_m}\right)$ of order $m$, and the rank-one tensor $x^m$:
$$f(x) \equiv \A x^m := \sum_{i_1, \cdots, i_m = 1}^n a_{i_1, \cdots, i_m}x_{i_1}
\cdots x_{i_m}.$$ The tensor $\A$ is called {\bf symmetric} as its
entries $a_{i_1, \cdots, i_m}$ are invariant under any permutation
of their indices. The tensor $\A$ is called {\bf positive definite
(semidefinite)} if $f(x) > 0$ ($f(x) \ge 0$) for all $x \in \Re^n, x
\not = 0$. When $m$ is even, the positive definiteness of such a
homogeneous polynomial form $f(x)$ plays an important role in the
stability study of nonlinear autonomous systems via Liapunov's
direct method in {\bf Automatic Control} \cite{NQW}. For $n$ big and
$m \ge 4$, this issue is a hard problem in mathematics.

In 2005, Qi \cite{Qi} defined eigenvalues and eigenvectors of a real
symmetric tensor, and explored their practical applications in
determining positive definiteness of an even degree multivariate
form.

By the tensor product, $\A x^{m-1}$ for a vector $x \in \Re^n$
denotes a vector in $\Re^n$, whose $i$th component is
$$\left(\A x^{m-1}\right)_i \equiv \sum_{i_2, \cdots, i_m = 1}^n a_{i, i_2, \cdots, i_m}x_{i_2} \cdots x_{i_m}.$$

We call a number $\lambda \in C$ an {\bf eigenvalue} of $\A$ if it
and a nonzero vector $x \in C^n$ are solutions of the following
homogeneous polynomial equation:
\begin{equation} \label{eig}
\left(\A x^{m-1}\right)_i = \lambda x_i^{m-1}, \ \ \forall \ i=1,
\cdots, n.
\end{equation}
and call the solution $x$ an {\bf eigenvector} of $\A$ associated
with the eigenvalue $\lambda$. We call an eigenvalue of $\A$ an {\bf
H-eigenvalue} of $\A$ if it has a real eigenvector $x$. An
eigenvalue which is not an H-eigenvalue is called an {\bf
N-eigenvalue}.   A real eigenvector associated with an H-eigenvalue
is called an {\bf H-eigenvector}.

The {\bf resultant} of (\ref{eig}) is a one-dimensional polynomial
of $\lambda$.   We call it the {\bf characteristic polynomial} of
$\A$.

We have the following conclusions on eigenvalues of an $m$th order
$n$-dimensional symmetric tensor $\A$:

\begin{theorem} ({\bf Qi 2005}) \label{t1}

(a). A number $\lambda \in C$ is an eigenvalue of $\A$ if and only
if it is a root of the characteristic polynomial $\phi$.

(b). The number of eigenvalues of $\A$ is $d=n(m-1)^{n-1}$.  Their
product is equal to det$(\A)$, the resultant of $\A x^{m-1} = 0$.

(c). The sum of all the eigenvalues of $A$ is
$$(m-1)^{n-1}{\rm tr}(\A),$$
where tr$(\A)$ denotes the sum of the diagonal elements of $\A$.

(d). If $m$ is even, then $\A$ always has H-eigenvalues.   $\A$ is
positive definite (positive semidefinite) if and only if all of its
H-eigenvalues are positive (nonnegative).

(e). The eigenvalues of $\A$ lie in the following $n$ disks:

$$ | \lambda - a_{i, i, \cdots, i} | \le \sum \{ | a_{i, i_2,
\cdots, i_m} | : i_2, \cdots, i_m = 1, \cdots, n, \{ i_2, \cdots,
i_m \} \not = \{ i, \cdots, i \} \}, $$ for $i = 1, \cdots, n$.
\end{theorem}

The resultant of $\A x^{m-1} = 0$, denoted as det$(\A)$ in this
theorem, is called the symmetric hyperdeterminant of $\A$.  When $m
\ge 3$, it is different from the hyperdeterminant,  introduced by
Cayley \cite{Ca} in 1845, and studied by Gelfand, Kapranov and
Zelevinsky \cite{GKZ} in 1994.  When $m=2$, both the
hyperdeterminant and the symmetric hyperdeterminant are the same as
the determinant in the classical sense.

This theorem extended the classic theory of matrix eigenvalue to
tensors.   However, eigenvalues and H-eigenvalues are not invariant
under orthogonal transformation.   This is needed in mechanics and
physics.    Hence, in \cite{Qi}, another kind of eigenvalues of
tensors was defined.

Suppose that $\A$ is an $m$th order $n$-dimensional symmetric
tensor. We say a complex number $\lambda$ is an {\bf E-eigenvalue}
of $\A$ if there exists a complex vector $x$ such that
\begin{equation} \label{eig2}
\left\{{\A x^{m-1} = \lambda x, \atop x^Tx = 1.}\right.
\end{equation}
In this case, we say that $x$ is an E-eigenvector of the tensor $\A$
associated with the E-eigenvalue $\lambda$.   If an E-eigenvalue has
a real E-eigenvector, then we call it a {\bf Z-eigenvalue} and call
the real E-eigenvector a {\bf Z-eigenvector}.

When $m$ is even, the {\bf resultant} of
$$\A x^{m-1} - \lambda (x^Tx)^{m-2 \over 2}x = 0$$
is a one dimensional polynomial of $\lambda$ and is called the {\bf
E-characteristic polynomial} of $\A$.  We say that $\A$ is regular
if the following system has no nonzero complex solutions:
$$\left\{{\A x^{m-1} = 0, \atop x^Tx = 0.}\right.$$

Let $P = (p_{ij})$ be an $n \times n$ real matrix.   Define $\B =
P^m\A$ as another $m$th order $n$-dimensional tensor with entries
$$b_{i_1, i_2, \cdots, i_m} = \sum_{j_1, j_2, \cdots, j_m = 1}^n p_{i_1j_1}p_{i_2j_2}\cdots p_{i_mj_m}a_{j_1, j_2, \cdots, j_m}.$$
If $P$ is an orthogonal matrix, then we say that $\A$ and $\B$ are
{\bf orthogonally similar}.

\begin{theorem} ({\bf Qi 2005}) \label{t2}

We have the following conclusions on E-eigenvalues of an $m$th order
$n$-dimensional symmetric tensor $\A$:

(a). When $\A$ is regular, a complex number is an E-eigenvalue of
$\A$ if and only if it is a root of its E-characteristic polynomial.

(b). Z-eigenvalues always exist.    An even order symmetric tensor
is positive definite if and only if all of its Z-eigenvalues are
positive.

(c). If $\A$ and $\B$ are orthogonally similar, then they have the
same E-eigenvalues and Z-eigenvalues.

(d). If $\lambda$ is the Z-eigenvalue of $\A$ with the largest
absolute value and $x$ is a Z-eigenvector associated with it, then
$\lambda x^m$ is the best rank-one approximation of $\A$, i.e.,
$$\| \A - \lambda x^m \|_F = \sqrt{ \| \A \|_F^2 - \lambda^2} = \min \{ \| \A - \alpha u^m \|_F
: \alpha  \in \Re, u \in \Re^n, \| u \|_2 = 1 \},$$ where $\| \cdot
\|_F$ is the Frobenius norm.
\end{theorem}

Theorem \ref{t2} (d) indicates that Z-eigenvalues play an important
role in the best rank-one approximation.  This also implies that
Z-eigenvalues are significant in practice.  The best rank-one
approximation of higher order tensors has extensive engineering and
statistical applications, such as {\bf Statistical Data Analysis}
\cite{DDV, KR, ZG}.

Eigenvalues, H-eigenvalues, E-eigenvalues and Z-eigenvalues were
extended to nonsymmetric tensors in \cite{Qi2}.  When $m$ is odd,
the E-characteristic polynomial was defined in \cite{Qi2} as the
{\bf resultant} of
$$\left\{{\A x^{m-1} - \lambda
x_0^{m-2}x, \atop x^\top x - x_0^2.}\right.$$

Independently, Lim also defined eigenvalues for tensors in
\cite{Li}. Lim defined eigenvalues for general real tensors in the
real field.     The $l^2$ eigenvalues of tensors defined by Lim are
Z-eigenvalues of Qi \cite{Qi}, while the $l^k$ eigenvalues of
tensors defined by Lim are H-eigenvalues in Qi \cite{Qi}. Notably,
Lim proposed a multilinear generalization of the Perron-Frobenius
theorem based upon the notion of $l^k$ eigenvalues (H-eigenvalues)
of tensors.

The calculation of eigenvalues is NP-hard with respect to $n$ when
$m \ge 3$.   In the case of Z-eigenvalues, when $n$ is small, we may
solve polynomial system (\ref{eig2}) by elimination methods, such as
the Gr\"obner method. When $n = 2$ or $3$, we may eliminate
$\lambda$ in $\A x^{m-1} = \lambda x$ first.  See \cite{QWW}.   When
$n$ is not small, we may use the power method. This method only can
get one eigenvalue and cannot guarantee it is the largest or the
smallest eigenvalue.   See \cite{KM}.

\medskip

{\bf B. Eigenvalues of Nonnegative Tensors, Higher-Order Markov
Chains and Spectral Hypergraph Theory}

This group include four research topics.

\smallskip

{\bf (1)  The Largest H-Eigenvalue of a Nonnegative Tensor.} We call
$\rho(\mathcal{A})$ the {\bf spectral radius} of tensor
$\mathcal{A}$ if
\begin{displaymath}
\rho(\mathcal{A})=\max\{|\lambda|:\, \text{$\lambda$ is an
eigenvalue of $\mathcal{A}$}\},
\end{displaymath}
where $|\lambda|$ denotes the modulus of $\lambda$.  A tensor
$\mathcal{A} = (a_{i_1\cdots i_m})$ is called {\bf nonnegative}, if
$a_{i_1\cdots i_m} \ge 0$, for all $i_1,\ldots,i_m = 1, \cdots, n.$
In the study of the Perron-Frobenius theorem for nonnegative tensors
and related algorithms, several classes of (nonnegative) tensors
were generalized to (nonnegative) tensors, some new classes of
(nonnegative) tensors were introduced. Among these, Chang, Pearson
and Zhang \cite{CPZ1} extended irreducible matrices to irreducible
tensors; Pearson \cite{Pe} extended essentially positive matrices to
essentially positive tensors;  Chang, Pearson and Zhang \cite{CPZ2}
extended primitive matrices to primitive tensors; Friedland, Gaubert
and Han \cite{FGH} introduced weakly irreducible tensors and wekly
primitive tensors; Zhang and Qi \cite{ZQ} introduced weakly positive
tensors; Hu, Huang and Qi \cite{HHQ} introduced strictly nonnegative
tensors and described the relations among corresponding nonnegative
tensors.

Now we may state the Perron-Frobenius theorem for nonnegative
tensors as follows.

\begin{theorem} {\bf (The Perron-Frobenius Theorem for Nonnegative Tensors)}
If $\A$ is a nonnegative tensor of order $m$ and dimension $n$, then
$\rho(\A)$ is an eigenvalue of $\A$ with a nonnegative eigenvector
$x \in \Re^n_+$.  (Yang and Yang \cite{YY1})

\medskip

If furthermore $\A$ is strictly nonnegative, then $\rho(\A)
> 0$. (Hu, Huang and Qi \cite{HHQ})

\medskip

If furthermore $\A$ is weakly irreducible, then $\rho(\A)$ is an
eigenvalue of $\A$ with a positive eigenvector $x \in \Re^n_{++}$.
(Friedland, Gaubert and Han \cite{FGH})

\medskip
Suppose that furthermore $\A$ is irreducible. If $\lambda$ is an
eigenvalue with a nonnegative eigenvector, then $\lambda=\rho(\A)$.
(Chang, Pearson and Zhang \cite{CPZ1})

\medskip
In this case, if there are $k$ distinct eigenvalues of modulus
$\rho(\A)$, then the eigenvalues are $\rho(\A)e^{i2\pi j / k}$,
where $j = 0, \cdots, k-1$.   The number $k$ is called the {\bf
cyclic index} of $\A$.   (Yang and Yang \cite{YY1})

\medskip

If moreover $\A$ is primitive, then its cyclic number is 1. (Chang,
Pearson and Zhang \cite{CPZ2})

\medskip

If $\A$ is essentially positive or even-order irreducible, then the
unique positive eigenvalue is real geometrically simple.  (Pearson
2010) (Yang and Yang \cite{YY})

\end{theorem}

Chang, Pearson and Zhang \cite{CPZ1} extended the well-known Collatz
minimax theorem for irreducible nonnegative matrices to irreducible
nonnegative tensors.  Based upon this, Ng, Qi and Zhou proposed an
algorithm for calculating $\rho(\A)$ in \cite{NQZ}.   For a weakly
irreducible nonnegative tensor $\A$, this algorithm generates two
sequence $\{ \underline{\lambda}_k \}$ and $\{ \bar{\lambda}_k \}$
such that
\[
\underline{\lambda}_1 \le \underline{\lambda}_2 \le \cdots \le
\lambda_0 = \rho(\A) \le \cdots \le  \bar{\lambda}_2 \le
\bar{\lambda}_1.
\]
The convergence of the NQZ algorithm was established by Chang,
Pearson and Zhang \cite{CPZ2} under primitivity, and by Friedland,
Gaubert and Han \cite{FGH} under weak primitivity.   The linear
convergence of the NQZ algorithm was established by Zhang and Qi
\cite{ZQ}.

Further developments in this topic can be found in \cite{Ch, LZI,
SQ, SQ1, YY3, YYL, Zh, ZQL, ZQX, ZCQ, ZCTW, ZQW}.

\smallskip

{\bf (2)  Higher-Order Markov Chains and and Z-eigenvalues of
Nonnegative Tensors.}  Ng, Qi and Zhou put higher-order Markov
chains as an application in their paper.   A higher-order Markov
chain model is used to fit the observed data through the calculation
of higher-order transition probabilities:
\begin{equation}\label{eq1}
0 \le p_{k_1k_2\cdots k_m} = {\rm Prob}(X_t = k_{1} \ | \
X_{t-1}=k_2, \ldots, X_{t-m+1}= k_m) \le 1
\end{equation}
where
\begin{equation} \label{sum1}
\displaystyle \sum_{k_1=1}^n p_{k_1,k_2,\cdots,k_m} = 1.
\end{equation}

Thus, we have an $m$-order $n$-dimensional nonnegative tensor ${\cal
P}$ consisting of $n^m$ entries in between 0 and 1:
$$
{\cal P} = (p_{k_1 k_2 ... k_m}), \quad 1\le k_1, k_2, ..., k_m \le
n.
$$
Let the probability distribution at time $t$ be $x^{(t)} \in S_n$,
where $S_n = \{ x \in \Re^n : x \ge 0, \sum_{j=1}^n x_j = 1 \}$.
Then we have
\begin{equation} \label{mark}
x^{(t)} = {\cal P}x^{(t-1)}\cdots x^{(t-m+1)} : =
\left(\sum_{k_2\cdots k_m=1}^n p_{ik_2\cdots
k_m}x^{(t-1)}_{k_2}\cdots x^{(t-m+1)}_{k_m}\right)_{i=1}^n \in S_n
\end{equation}
for $t = m+1, \cdots$. Assume that
\begin{equation} \label{station}
\lim_{t\to\infty} x^{(t)} = x^*.
\end{equation}
Then we have $x^* \in S_n$ and we may call $x$ as the stationary
probability distribution of the higher-order Markov chain.   We see
$x^*$ satisfies the following tensor equation:
\begin{equation} \label{fix}
{\cal P}x^{m-1} = x \in S_n.
\end{equation}

Li and Ng \cite{LN} showed that if ${\cal P}$ is irreducible, then a
positive solution of (\ref{fix}) exists.  Then also give sufficient
conditions for uniqueness of the solution of the {\bf fixed point
problem} (\ref{fix}).  Under the same conditions, Li and Ng
\cite{LN1} established linear convergence of the power method
\begin{equation} \label{power}
y^{(k+1)} = {\cal P}\left(y^{(k)}\right)^{m-1}
\end{equation}
Li and Ng \cite{LN1} established linear convergence of the power
method (\ref{power}), but has not proved the convergence of the
higher order Markov chain (\ref{mark}), i.e., (\ref{station}).   For
second-order Markov chain, under some conditions, Hu and Qi
\cite{HQ1} established the convergence (\ref{station}).

The fixed point problem (\ref{fix}) is different from the eigenvalue
problem (\ref{eig}). It is somewhat connected with the Z-eigenvalue
problem of a special class of nonnegative tensors.  Unlike the
largest H-eigenvalue of a nonnegative tensor, the largest
Z-eigenvalues of a nonnegative tensor are not unique in general
\cite{CPZ1}.  The recent paper of Chang, Pearson and Zhang
\cite{CPZ3} reveals some similarities as well as differences between
the Z-eigenvalues and H-eigenvalues of a nonnegative tensor.

\smallskip

{\bf (3) Spectral Hypergraph Theory.} Lim \cite{Li} pointed out that
a potential application of the largest eigenvalue of a nonnegative
tensor is on {\bf hypergraphs}. A graph can be described by its
adjacency matrix. The properties of eigenvalues of the adjacency
matrix of a graph are related to the properties of the graph. This
is the topic of spectral graph theory. A hypergraph can be described
by a (0, 1)- symmetric tensor, which is called its adjacency tensor.
Are the properties of the eigenvalues of the adjacency tensor of a
hypergraph related to the properties of the hypergraph?  Recently,
several papers appeared in spectral hypergraph theory \cite{BP, BP1,
CD, HQ, LQY1}.  It is interesting that in \cite{BP, BP1, CD},
H-eigenvalues are used, while in \cite{HQ, LQY1}, Z-eigenvalues are
used. This shows that there are more to explore in this topic.    We
expect a comprehensive spectral hypergraph theory will emerge
eventually.

\smallskip

{\bf (4)  The Relation with Nonnegative Tensor Factorization.}
Nonnegative tensor factorization is an important topic in the tensor
decomposition community \cite{CZPA, SH}.   What is the relation
between the spectral theory of nonnegative tensors with the
nonnegative tensor factorization topic?   Currently, it is still
blank in the literature.

\medskip

{\bf C. Eigenvalues of General Tensors: Theory}

This group include six research topics.

\smallskip

{\bf (1)  The Number of E-Eigenvalues, the Number of Real
Eigenvalues, E-Characteristic Polynomial.} In \cite{NQWW}, Ni, Qi,
Wang and Wang showed that the degree of the E-characteristic
polynomial of an even order tensor is bounded by
$$d(m, n) = {(m-1)^n - 1 \over m-2},$$
where $m$ is the order and $n$ is the dimension, $m$ is even and $m
\ge 4$.    The definition of E-eigenvalues in \cite{Qi} is not a
strict generalization of eigenvalues of a square matrix, as it
excludes the case that $x^\top x = 0$.    In \cite{CS}, Cartwright
and Sturmfels gave a precise definition by using equivalence
classes, and showed the number of the equivalent classes is
$$d(m, n) = {(m-1)^n - 1 \over m-2},$$
for any $m > 2$, echoing the result of \cite{NQWW}.

In \cite{CPZ}, Chang, Pearson and Zhang unified the definitions of
H-eigenvalues, Z-eigenvalues and D-eigenvalues, and showed that for
a real $n$-dimensional $m$th order symmetric tensor, there are at
least $n$ H-/Z-/D-eigenvalues.   Existence of real eigenvalues of
real tensors was further studied in \cite{Zh1}.

In \cite{LQZ}, Li, Qi and Zhang studied the properties of
E-characteristic polynomials.

\smallskip

{\bf (2)  Symmetric Hyperdeterminant or E-determinant.}  In
\cite{CD}, Cooper and Dutle gave explicit expression of symmetric
hyper-determinants of some special classes of hypergraphs for
arbitrary $n$.   This shows that for some special cases, symmetric
hyper-determinants are computable.    In \cite{LQZ}, Li, Qi and
Zhang showed that symmetric hyper-determinants are invariant under
orthogonal transformation. In \cite{HHLQ}, symmetric
hyper-determinants are called E-determinants.  By the definition,
det$(\A) = 0$ if and only if $\A x^{m-1} = 0$ has a nonzero
solution.   Hu, Huang, Ling and Qi \cite{HHLQ} showed that if
det$(\A) \not = 0$, then any nonhomogeneous polynomial system with
$\A$ as its leading coefficient tensor has a solution. This shows
that like determinants are closely related to the solution theory of
linear systems, hyper-determinants are closely related to the
solution theory of polynomial systems.  More research is needed in
this topic.

\smallskip

{\bf (3)  Singular Values, Symmetric Embedding and Eigenvalue
Inclusion.} Singular values of rectangular tensors were studied in
\cite{CQZ, CT, Li, YY3, Zh, ZCQ}.   Symmetric embedding, which links
singular values of rectangular values, and eigenvalues of square
tensors, was studied in \cite{CL, RV}.

In \cite{LL, LLK}, the authors give a number of eigenvalue inclusion
theorems, including a Taussky-type boundary result and a Brauer
eigenvalue inclusion theorem for tensors, and their applications.

\smallskip

{\bf (4)  Waring Decomposition Related Problems.} A Waring
decomposition of a (homogeneous) polynomial $f$ is a minimal sum of
powers of linear forms expressing $f$. Under certain conditions,
such a decomposition is unique. In \cite{OO}, Oeding and Ottaviani
discussed some algorithms to compute the Waring decomposition, which
are linked to the equations of certain secant varieties and to
eigenvectors of tensors. In particular they explicitly decomposed a
cubic polynomial in three variables as the sum of five cubes
(Sylvester Pentahedral Theorem).

\smallskip

{\bf (5) Third-Order Tensors as Linear Operators.} In \cite{Br, GER,
KBHH}, special operations were introduced to regard third-order
tensors as linear operations.   With more complicated manipulations,
this approach may also be extended to higher-order tensors.  The
cost of this approach is that the original tensor operations are
given up.

\smallskip

{\bf (6) Successive Symmetric Best Rank-One approximation, the Best
Rank-One Approximation Ratio, and the Symmetric Rank Conjecture.}
Denote the space of all real $m$th-order $n$ dimensional symmetric
space as $S_{m, n}$. Suppose that $\A \in S_{m, n}$. Let $\A^{(0)}
\equiv \A$.   In general, suppose that $\A^{(k)} \in S_{m, n}$. Let
$\lambda_k$ be the Z-eigenvalue of $\A^{(k)}$ with the largest
absolute value, and $x^{(k)}$ be the corresponding Z-eigenvector.
Then we define
$$\A^{(k+1)} = \A^{(k)} - \lambda_k \left(x^{(k)}\right)^m.$$
Wang and Qi \cite{WQ} proved that
\begin{equation} \label{succ}
\A = \sum_{k=0}^\infty \lambda_k \left(x^{(k)}\right)^m,
\end{equation}
and called this the successive symmetric best rank-one approximation
of $\A$.

Qi \cite{Qi3} studied further the convergence of (\ref{succ}).
Denote the largest absolute value of the Z-eigenvalues of $\A \in
S_{m, n}$ as $\rho_Z (\A)$.   When $m=2$, $\rho_Z (\A)$ is actually
the spectral radius, i.e., the 2-norm of $\A$.   In $S_{m, n}$,
$\rho_Z (\A)$ is still a norm of $\A$.  Denote the Frobenius norm of
$\A$ by $\| \A \|_F$.  Then
$$\App (S_{m, n}) := \min \left\{ {\rho_Z (\A) \over \| \A \|_F} : \A \in S_{m, n} \right\} > 0.$$
Qi \cite{Qi3} called it the {\bf best rank-one approximation ratio}
of $S_{m, n}$.  Actually, we have
$$0 < \App (S_{m, n}) < 1.$$
Suppose that $\A^{(k)} \in S_{m, n}$.   Let $\lambda_k$ be the
Z-eigenvalue of $\A^{(k)}$ with the largest absolute value, and
$x^{(k)}$ be the corresponding Z-eigenvector. Then we have
$$\left\| \A^{(k+1)} \right\|_F^2 \le \left\| \A^{(k)} \right\|_F^2 \left[ 1 - \App (S_{m,
n})^2\right].$$ This not only proves (\ref{succ}), but also gives
its convergence rate.

Let
$$\underline{\mu}_{m,n} = {1 \over \sqrt{n^{m-1}}}.$$

If $m=2k$ is even, then let $\A^{(m, n)} \in S_{m, n}$ and ${\bar
\mu}_{m, n}$ be defined by
\begin{equation} \label{evenA}
\A^{(m,n)} x^m = \left(x^\top x\right)^k
\end{equation}
and
$${\bar \mu}_{m, n} = { 1 \over \|\A^{(m,n)} \|}.$$

If $m=2k+1$ is odd, then let $\A^{(m, n)} \in S_{m, n}$ and ${\bar
\mu}_{m, n}$ be defined by
\begin{equation} \label{oddA}
\A^{(m,n)} x^m = \left(x^\top x\right)^k\left(\sum_{i=1}^n
x_i\right)
\end{equation}
and
$${\bar \mu}_{m, n} = { \sqrt{n} \over \|\A^{(m,n)} \|}.$$

The following upper and lower bounds for $\App (S_{m, n})$ were
established in \cite{Qi3}.

\begin{theorem} \label{thm4.4}
The value  $\underline{\mu}_{3,n}$ is a positive lower bound for
$\App (S_{3, n})$.

On the other hand, the value  ${\bar \mu}_{m, n}$ is an upper bound
for ${\App} (S_{m, n})$ for $m =2, 3, \cdots$. We have
\begin{equation} \label{mu1}
{1 \over \sqrt{n}} = \underline{\mu}_{2,n} = {\App} (S_{2, n}) =
{\bar \mu}_{2,n} = {1 \over \sqrt{n}},
\end{equation}
\begin{equation} \label{mu2}
{1 \over n} = \underline{\mu}_{3, n} \le {\App} (S_{m, n}) \le {\bar
\mu}_{3, n} = \sqrt{6 \over n + 5}
\end{equation}
and
\begin{equation} \label{mu3}
{\App} (S_{4, n}) \le {\bar \mu}_{4, n} = \sqrt{3 \over n^2 + 2n}.
\end{equation}
\end{theorem}

To improve this theorem, Qi \cite{Qi3} made a conjecture that
$\rho_Z (\A)$ is equal to
$$\sigma (\A) : = \max \left\{ \left|\A x^{(1)}\cdots x^{(m)}\right| :
\left(x^{(i)}\right)^\top x^{(i)} = 1, \forall i = 1, \cdots, m
\right\}.$$ This conjecture was proved by Zhang, Ling and Qi
\cite{ZLQ2}.  This improved Theorem \ref{thm4.4} that
$\underline{\mu}_{m,n}$ is a positive lower bound for $\App (S_{m,
n})$ for all $m$.

One referee of \cite{ZLQ2} pointed out the above result can be
regarded the first step to prove the symmetric rank conjecture of
Comon, Golub, Lim and Mourrain \cite{CGLM}.

Let $\A$ be an $m$th-order $n$-dimensional tensor.  We may have a
rank-one tensor decomposition of $\A$:
$$\A = \sum_{i=1}^r \alpha_i x^{(i, 1)}\cdots x^{(i, m)},$$
where $x^{(i, j)}$ are $n$-dimensional vectors, $\alpha_i$ are
numbers.   The minimum value of $r$ is called the rank of $\A$.

If $\A$ is symmetric, than we may have a symmetric rank-one tensor
decomposition of $\A$:
$$\A = \sum_{i=1}^r \alpha_i \left(x^{(i)}\right)^m,$$
where $x^{(i)}$ are $n$-dimensional vectors, $\alpha_i$ are numbers.
The minimum value of $r$ is called the symmetric rank of $\A$.

It is conjectured by Comon, Golub, Lim and Mourrain \cite{CGLM} that
for a symmetric tensor, its symmetric rank is the same as its rank.
We call this conjecture the {\bf symmetric rank conjecture}.

\medskip

{\bf D. Eigenvalues of General Tensors: Links with Optimization,
Numerical Analysis and Geometry}

This group include four research topics.

\smallskip

{\bf (1)  Bi-quadratical Optimization and Spherical Optimization.}
Bi-quadratical optimization and spherical optimization were studied
in \cite{BLQZ, HLZ, HH, LNQY, LZQ, LZ, So, WQZ, ZLQ, ZLQ2, ZQY}.
Bi-quadratical optimization is linked with the M-eigenvalue problem
\cite{DQH, HQD, LNQY, WQZ}, while spherical optimization is linked
with the Z-eigenvalue problem.   They are NP-hard in general.   The
approach in \cite{BLQZ, HLZ, LNQY, LZQ, LZ, So, ZLQ, ZLQ2} is to use
semi-definite programming to find a lower bound for the minimization
problem.   Symmetric spherical optimization was studied in
\cite{ZLQ2, ZQY}.  As stated in the last section, Zhang, Ling and Qi
\cite{ZLQ2} proved that for $\A \in S_{m, n}$,
$$\rho_Z (\A) \equiv \max \left\{ \left|\A x^m \right| :
x^\top x = 1  \right\}$$
$$= \sigma (\A) \equiv \max \left\{ \left|\A x^{(1)}\cdots x^{(m)}\right| :
\left(x^{(i)}\right)^\top x^{(i)} = 1, \forall i = 1, \cdots, m
\right\}.$$ This reveals some deep insight property of symmetric
spherical optimization.

\smallskip

{\bf (2)  Space Tensor Conic Programming.} In the study of higher
order diffusion tensor imaging, positive semi-definite tensor models
arise \cite{HHNQ, QYX, QYW}, where the tensor is positive definite
and of dimension three.  We call three dimensional tensors {\bf
space tensors}.  Space tensors appear in physics and mechanics. They
are real physical entities.   The positive semi-definite space
tensors of the same order form a convex cone.   The optimization
problem involving positive semi-definite space tensors is called the
{\bf space tensor conic programming} problem \cite{QY}.  In
\cite{LQY1}, based on the analysis of the semismoothness properties
of the maximum Z-eigenvalue function, a generalized Newton method is
proposed to solve the space tensor conic linear programming problem.
In \cite{HHQ1}, a sequential SDP method is proposed to solve the
space tensor conic linear programming problem.

The problem to identify a space tensor is positive (semi-)definite
or not is equivalent to finding the minimum Z-eigenvalue of an
$m$-th order tensor with dimension three.  As we discussed in the
basic theory section, this problem can be solved by elimination
methods. It is not NP-hard.   This approach requests to solving a
one variable polynomial equation of degree generically $N = {1 \over
2} (m+1)(m+2)$.  When $m$ is big, the solution procedure is
unstable. Hence, new algorithms are needed and this research topic
is interesting.

\smallskip

{\bf (3) Characterizing the Limiting Behaviour of Newton's Method.}
Dupont and Scott \cite{DS} studied the angular orientations of
convergent iterates generated by Newton's method in multiple space
dimensions. They showed that the Newton iteration can be interpreted
as a fixed-point algorithm for solving a tenor eigenproblem.  They
gave an extensive computational analysis of this tensor eigenproblem
in two dimensions. In a large fraction of cases, the tensor
eigenproblem has a discrete number of solutions to which the Newton
directions converge quickly, but there is also a large fraction of
cases in which the behavior is more complicated.  Dupont and Scott
contrasted the angular orientations of iterates generated by
Newton's method with the corresponding directions of the continuous
Newton algorithm.

\smallskip

{\bf (4) Geometry related Problems.}   Qi \cite{Qi1} studied the
orthogonal classification problem of real hypersurfaces given by an
equation of the form
$$S = \{ x \in \Re^n: f(x) = \A x^m = c \},$$
where $\A$ is a symmetric tensor, through the rank, the
Z-eigenvalues and the asymptotic directions.

Balan \cite{Ba, Ba1, Ba2, BP} showed that there is a close
relationship between Finsler metrics of Berwald-Moor type, and the
spectral theory of tensors. Balan used the results of \cite{Qi, Qi1}
to analyze spectral properties of m-root models, geometric relevance
of spectral equations, degeneracy sets, asymptotic rays, base
indices and based rank one approximation of various Berwald-Moor
tensors, Chernov tensors and Bogoslovski tensors.    A further
investigation on the geometrical meanings of eigenvalues of higher
order tensors is needed.

\medskip

{\bf E. Eigenvalues of General Tensors: Links with Image Science,
Solid Mechanics and Quantum Physics}

This group include four research topics.

\smallskip

{\bf (1)  Higher Order Diffusion Tensor Imaging.}  Diffusion tensor
imaging (DTI) is the most popular magnetic resonance imaging model
\cite{BJ}.  In the DTI model, eigenvalues of the diffusion tensors
play a fundamental role as they are the main invariants of these
tensors, and thus can be used in the measures of DTI \cite{BJ}.
However, DTI is known to have a limited capability in resolving
multiple fibre orientations within one voxel.  This is mainly
because the probability density function for random spin
displacement is non-Gaussian in the confining environment of
biological tissues and, thus, the modeling of self-diffusion by a
second order tensor breaks down.   Hence, researchers presented
various {\bf Higher Order Diffusion Tensor Imaging} models to
overcome this problem \cite{JHRLK, LJRH, OM, Tu, TRWMBW}.  For the
measures of these higher order diffusion tensor imaging models, the
spectral theory of tensors naturally play an important roles
\cite{BV, CHCHQW, HQW, HHNQ, HCQW, HCQW1, QHW, QWW1, QYX, QYW,
ZLQW}.  In particular, Qi, Wang and Wu \cite{QWW1} proposed
D-eigenvalues to study the diffusion kurtosis model \cite{JHRLK,
LJRH}; Bloy and Verma \cite{BV} studied the underlying fiber
directions from the diffusion orientation distribution function by
using Z-eigenvalues; Qi, Yu and Wu \cite{QYW}, Hu, Huang, Ni and Qi
\cite{HHNQ}, and Qi, Yu and Xu \cite{QYX} studied positive
semi-definite higher order diffusion tensor imaging models.

\smallskip

{\bf (2)  Image Authenticity Verification Problem.} In \cite{ZZP},
Zhang, Zhou and Peng generalized the definition of D-eigenvalues,
and introduced the gradient skewness tensor which involves a
three-order tensor derived from the skewness statistic of gradient
images. They found out that the skewness value of oriented gradients
of an image can measure the directional characteristic of
illumination, and the local illumination detection problem for an
image can be abstracted as solving the largest D-eigenvalue of
gradient skewness tensors. Their method presented excellent results
in a class of image authenticity verification problem, which is to
distinguish real and flat objects in a photograph.

\smallskip

{\bf (3)  Elasticity Tensors in Solid Mechanics.} In solid
mechanics, the elasticity tensor $\A = (a_{ijkl})$ is partially
symmetric in the sense that  for any $i, j, k, l$, we have $a_{ijkl}
= a_{kjil} = a_{ilkj}$.  We say that they are strongly elliptic if
and only if
$$f(x, y) \equiv \A xyxy \equiv \sum\limits_{i,j,k,l=1}^n a_{ijkl}x_iy_jx_ky_l > 0,$$
for all unit vectors $x, y\in \Re^n, n=2~{\rm or}~ 3$. For an
isotropic material, some inequalities have been established to judge
the strong ellipticity \cite{KS, WA}.  In \cite{HDQ, QDH}, Dai, Han
and Qi studied conditions for strong ellipticity and introduced {\bf
M-eigenvalues} for the ellipticity tensor $\A$.  They showed that
M-eigenvalues always exist and the strong ellipticity condition
holds if and only if the smallest M-eigenvalue of the elasticity
tensor is positive.

\smallskip

{\bf (4)  Quantum Entanglement Problem.} The entanglement problem is
to determine whether a quantum state is separable or inseparable
(entangled), or to check whether an $mn \times mn$ symmetric matrix
$A\succeq 0$ can be decomposed as a convex combination of tensor
products of $n$ and $m$ dimensional vectors. It has fundamental
importance in quantum science and has attracted much attention since
the pioneer work of Einstein, Podolsky and Rosen \cite{EPR} and
Schr\"odinger \cite{Sc}.  The typical optimization problem in
quantum entanglement problem is like the following \cite{HS}:
$$\max \left\{ \left|\A x^{(1)}\cdots x^{(m)}\right| : \left| x^{(i)}\right|
= 1, \ {\rm for } \ i = 1, \cdots, m \right\},$$ where $\A$ is an
$m$th-order $(n_1, \cdots, n_m)$-dimensional complex tensor,
$x^{(i)} \in C^{n_i}$ for $i = 1, \cdots, m$.  The singular value
form of this optimization problem is
$$\A x^{(1)}\cdots x^{(i-1)}x^{(i+1)}\cdots x^{(m)} = \lambda \bar
x^{(i)}, \ \left| x^{(i)} \right| = 1, \ {\rm for} \ i = 1, \cdots,
m.$$ The optimization problem and the singular value problem involve
complex and conjugate numbers.   They are different but have some
similarities with what we studied before.  With their strong physics
background, this topic will be a new fertile land of our research.

\medskip

The International Conference on The Spectral Theory of Tensors will
be held at Chern Research Institute, Nankai University, Tianjin,
China, during May 30 - June 2, 2012. Its website is:

http://www.nim.nankai.edu.cn/activites/conferences/hy20120530/index.htm

Two special issues on the Spectral Theory of Tensors are being
edited \cite{LNQ, YZZZ}.

\end{document}